\documentclass{amsart} 
\usepackage{amsmath,amscd} 
\usepackage{latexsym,amsmath,amssymb,mathrsfs,setspace,enumerate}
\usepackage[T1]{fontenc}
\usepackage[utf8]{inputenc}
\usepackage{lmodern} 
\usepackage{amssymb}
\usepackage{tikz-cd}
\usepackage{setspace}
\usepackage{graphics}
\usepackage{graphicx}
\usepackage[utf8]{inputenc}
\usepackage[english]{babel}
\usepackage{fancyhdr}
\usepackage{tikz}
\usetikzlibrary{decorations.markings,positioning}
\theoremstyle{plain}

\newtheorem{prop}{Proposition} 

\newtheorem{exam}{Example}

\newtheorem{defn}{Definition}

\tikzcdset{%
	triple line/.code={\tikzset{%
			double distance = 3pt, 
			double=\pgfkeysvalueof{/tikz/commutative diagrams/background color}}},
	quadruple line/.code={\tikzset{%
			double distance = 5.3pt, 
			double=\pgfkeysvalueof{/tikz/commutative diagrams/background color}}},
	Rrightarrow/.code={\tikzcdset{triple line}\pgfsetarrows{tikzcd implies cap-tikzcd implies}},
	RRightarrow/.code={\tikzcdset{quadruple line}\pgfsetarrows{tikzcd implies cap-tikzcd implies}}
}

\title{Chain bundles and Category of Chains}
\author {P. G. Romeo and Riya Jose}
\address{Dept. of Mathematics, Cochin University of Science and Technology, Kochi, Kerala, INDIA.}
\email{$romeo_-parackal@yahoo.com,\, riyajosemarangattu@gmail.com $}
\subjclass{20M10}
\keywords { Category, Subobject, Chainbundle, Cochain bundle, Chain bundle map, Factorization, Category of chains.}
\thanks{} 
\date{}

\begin{document}
	
	\begin{abstract}
		In \cite{PG}, we introduced  the category of chain buldles with examples and established its significance. Here we shall describe certain categories which we call the category of chains in the chain bundle category and discuss some interesting categorical properties of this chain categories. 
	\end{abstract}
	\maketitle
	
	Category theory, in general is a mathematical theory of structures and system of structures which has come to occupy a central position in  present day mathematics. Eilenberg and Mac Lane provide a purely abstract definition of a category in 1945 (cf. \cite{Samuel} and \cite{Mac}). 
	Accordingly a category $\mathcal{C}$ can be described as a set \textbf{Ob} whose members are objects of $\mathcal{C}$ and for every pair $X, Y$ of objects of $\mathcal{C}$ there is $Hom(X,Y)$ called morphisms from $X$ to $Y$ both satisfying certain specified conditions. 
	There are many situations where one have to deal with chains (sequences) of objects and hence it is natural to consider categories whose objects are chain of objects in some category $\mathcal{C}$ and morphisms are appropriate morphisms between such chains. In order to discuss such situations, in \cite{PG} we introduced the category of chain bundles as a category whose objects are sequence of objects in a category $\mathcal{C}$ with zero such that for any two objects $X$ and $Y$ any subset of $Hom(X,Y)$ constitute the morphisms between $X$ and $Y$. A category whose objects are chain bundles and morphisms appropriate chain bundle maps is termed as a category of chain bundles and will be denoted as $\mathfrak{CB_\mathcal{C}}$.\\
	
	In this paper we consider the category of chain bundles (cf. \cite{PG}) and descibe the chain bundle categories from directed graph and augmented simplicial categories. Further we discuss certain categorical properties of the chain bundle categories with examples and then describe the category of chains of chain bundle categories.
	\section{Preliminaries}
	In the following we briefly recall some basic notions such as categories, category with subobjects, factorization of morphisms and the like which are needed in the sequel.
	\begin{defn} A category $\mathcal{C}$ consists of
		\begin{enumerate}
			\item a class $ \nu\mathcal{C} $ called the class of vertices or objects
			\item a class $\mathcal{C}$ of disjoint sets $\mathcal{C}(a,b)$, one for each pair $(a,b) \in \nu\mathcal{C} \times  \nu\mathcal{C}  $, an element $f \in \mathcal{C}(a,b)$ called a morphism (arrow) from $a$ to $b$, where $a = dom\,f$ and $b = cod \,f$
			\item For $a,b,c \in  \nu\mathcal{C}$, a map 
			$$ \circ : \mathcal{C}(a,b) \times \mathcal{C}(b,c) \rightarrow \mathcal{C}(a,c)$$
			$$ (f,g) \rightarrow f\circ g \quad \text{is the composition of morphisms in \, } \mathcal{C}$$
			\item For each $ a \in \nu\mathcal{C}$, a unique $1_a \in \mathcal{C}(a,a)$ is the identity morphism on $a$
		\end{enumerate}
		These must satisfy the following axioms:
		\begin{enumerate}
			\item[Cat 1] The compostion is associative: for $f \in \mathcal{C}(a,b)$, $g \in \mathcal{C}(b,c)$ and $h \in \mathcal{C}(c,d)$; ie.,
			$$ f\circ(g\circ h) = (f\circ g)\circ h.$$
			\item[Cat 2] For each $ a \in \nu\mathcal{C},\, f \in \mathcal{C}(a,b),\, g \in \mathcal{C}(c,a),\quad 
			1_a \circ f =f \,\,\text{and}\,\,  g \circ 1_a = g.$		
		\end{enumerate}
	\end{defn}
	
	$\nu\mathcal{C}$ may identify as a subclass of $\mathcal{C}$ and with this identification, categories may regard 
	in terms of morphisms(arrows) alone. The category $\mathcal{C}$ is said to be small if the class $\mathcal{C}$ is 
	a set. For any category $\mathcal{C}$ an opposite category denoted as  $\mathcal{C}^{op} $ is the category with 
	$$ \nu \mathcal{C}^{op} = \nu \mathcal{C}, \mathcal{C}^{op}(a,b) = \mathcal{C}(b,a)\,\text{for all}\, a,
	b \in \nu \mathcal{C}$$
	ie., for arrows $f^{op}$ there is a one-one correspondance $f\mapsto f^{op}$ with arrow $f$ of $\mathcal{C}$ and the composite $f^{op}g^{op}=(gf)^{op}$ is defined in $\mathcal{C}^{op}$ exactly when $gf$ is defined in $\mathcal{C}$ (see cf. \cite{Mac}). 
	
	\begin{defn}
		A  covariant functor $F :\mathcal{C} \rightarrow \mathcal{D} $  from a category $\mathcal{C}$ to a category $\mathcal{D}$ consists of a \textit{vertex map} $ \nu F : \nu \mathcal{C} \rightarrow \nu \mathcal{D}$ which assigns to each $ a \in \nu \mathcal{C}$, a vertex $\nu F(a) \in \nu \mathcal{D}$ and a \textit{morphism map} $F$ which assigns to each morphism $f: a \rightarrow b $ in $ \mathcal{C}$, a morphism $F(f) : F(a) \rightarrow  F(b) \in \mathcal{D}$ such that 
		\begin{enumerate}
			\item[(Fn. 1)]  $F(1_a) = 1_{F(a)} \  \forall a \in \nu \mathcal{C}$
			\item[(Fn .2)] $F(f)F(g) = F(fg)$ for all morphisms $f, g \in \mathcal{C}$ for which the composite $fg$ exists.
		\end{enumerate}
		$F$ is a contravariant functor if $\nu F$ is as above and the morphism map assigns to each $f: a \rightarrow b $ in 
		$\mathcal{C}$, a morphism $F(f) : F(b) \rightarrow F(a) \in \mathcal{D}$ such that they satisfy axiom 
		(Fn. 1) and (Fn .2)$^*$ ie., $F(g)F(f) = F(fg) $ for all morphisms $f, g \in \mathcal{C}$ for which the 
		composite $fg$ exists.
	\end{defn}
	
	\begin{defn}
		A category $\mathcal{D}$ is a subcategory of a category $\mathcal{C}$ if the class $\mathcal{D}$ is a subclass of $\mathcal{C}$ and the composition in $\mathcal{D}$ is the restriction of the composition in $\mathcal{C}$ to $\mathcal{D}$. In this case, the inclusion $\mathcal{D} \subseteq \mathcal{C}$ preserves composition and identities and so represents a functor of $\mathcal{D}$ to $\mathcal{C}$. 
	\end{defn}
	
	\begin{defn} 
		A morphism $f$ in a category $\mathcal{C}$ is a monomorphism if 
		$$ gf = hf \Rightarrow g = h \ \forall g,h \in \mathcal{C} $$
		and a morphism $f$ in  $\mathcal{C}$ is called a split monomorphism if it has a right inverse. Every morphism in a concrete category whose underlying function is an injection is a monomorphism. A morphism $f$ in a category 
		$\mathcal{C}$ is an \textit{epimorphism} if 
		$$ fg = fh \Rightarrow g = h \ \forall g,h \in \mathcal{C} $$
		and a morphism $f$ in  $\mathcal{C}$ is called a split epimorphism if it has a left inverse. Every morphism in a concrete category whose underlying function is an surjection is an epimorphism.
	\end{defn}
	
	Let $\mathbb{M} \mathcal{C}$ denote the class of all monomorphisms in $\mathcal{C}$. On  $\mathbb{M} \mathcal{C}$, define the relation $$f \preceq g \Leftrightarrow f = hg \quad \text{for some } h \in \mathcal{C}$$
	$\preceq$ is a quasi- order and  
	\begin{equation}
		\sim = \preceq \cap \preceq^{-1}
	\end{equation} is an equivalence relation on $\mathbb{M} \mathcal{C}$.
	
	\begin{defn}
		A preorder $\mathcal{P}$ is a category such that for any $p, p' \in  \mathcal{P}$, $\mathcal{P}(p,p')$ contains atmost one morphism. In this case, the relation $\subseteq$ on the class $\nu \mathcal{P}$ defined by 
		\begin{equation}
			p \subseteq p' \ \Leftrightarrow \ \mathcal{P}(p,p') \neq \emptyset
		\end{equation}
		is a quasiorder. When $\mathcal{P}$ is a preorder, $\nu \mathcal{P}$ will stand for the quasiordered class $(\nu \mathcal{P}, \subseteq)$. Conversely given a quasiorder $\leq$ on the class $X$, the subset $$ \mathcal{P} = \{(x,y) \in X \times X : x \leq y\}$$
		of $ X \times X $ is a preorder such that the quasiordered class $\nu \mathcal{P}$ defined above is order isomorphic with $(X,\leq)$. If the relation $\subseteq$ on $\mathcal{P}$  is antisymmetric then we shall say that $\mathcal{P}$ is a \textit{strict preorder}.
	\end{defn}
	
	\begin{defn} 
		Let $ \mathcal{C} $ be a category. A choice of subobjects in $ \mathcal{C} $ is a subcategory $ \mathcal{P} \subseteq \mathcal{C} $ satisfying the following:
		\begin{enumerate}
			\item[(a)] $\mathcal{P}$ is a strict preorder with $\nu \mathcal{P} = \nu \mathcal{C}$.
			\item[(b)] Every $f \in \mathcal{P} $ is a monomorphism in  $\mathcal{C}$.
			\item[(c)] If $f, g \in \mathcal{P}$ and if $f = hg$ for some $h \in \mathcal{C}$ then $h \in \mathcal{P}$.
		\end{enumerate}
		When $\mathcal{P}$ satisfies these conditions, the pair $ (\mathcal{C},\mathcal{P})$ is called a caegory with subobjects.
	\end{defn}
	When $\mathcal{C} $ has  subobjects, unless explicitly stated otherwise, $ \nu\mathcal{C} $ will denote the choice of subobjects in $ \mathcal{C} $. The partial order defined by equation $(2)$ is called the preorder of inclusions or subobject relation in $ \mathcal{C} $ and is denoted by $ \subseteq $. If $c,d \in \nu\mathcal{C}$ and $ c\subseteq d $ the unique morphism from $c$ to $d$ is the inclusion $j_c^d: c \rightarrow d $. Any monomorphism $f$ equivalent to an inclusion (with respect to the equivalence relation $\sim$ defined by equation $(1)$) is called an \textit{embedding}.
	
	\begin{exam}
		In categories $Set, Grp, Vct_K, Mod_R $ the relation on objects induced  by the usual set inclusion is a subobject relation.
	\end{exam}
	
	\begin{defn}
		A morphism $f$ in a category $\mathcal{C}$ with subobjects is said to have \textit{factorization} if $f$ can be expressed as  $ f = pm$ where $p$ is an epimorphism and $m$ is an embedding.
	\end{defn}
	Factorization of a morphism need not be unique. Every morphism $f$ with factorization has atleast one factorization of the form $f = qj$  where $q$ is an epimorphism and $j$ is an inclusion. Such factorizations are called \textit{canonical factorization}. A category $\mathcal{C}$ is \textit{category with factorization} if $\mathcal{C}$ has subobjects and if every morphism in $\mathcal{C}$ has factorization. The category has \textit{unique factorization property} if every morphism in $\mathcal{C}$ has unique canonical factorization.
	
	\begin{exam}
		If $f : X \rightarrow Y$ is a mapping of sets and $f(X) = Im f$ the $f(X) \subseteq  Y$ and we can write $f = f^0j_{f(X)}^Y$. Here $f^0$  denote the mapping of $X$ onto $f(X)$ determined by $f$. Since surjective mappings are epimorphisms in $Set$, this gives a canonical factorization of $f$ in $Set$ which is clearly unique. ie., $Set$ is a category with unique factorization.
	\end{exam}
	
	\begin{prop}
		Let $\mathcal{C}$ be category with factorization. Suppose that the morphism $f \in \mathcal{C}$ has the following property:
		\begin{itemize}
			\item[(Im)] $f$ has a canonical factorization $f = xj$ such that for any canonical factorization $f = yj'$ of $f$, there is an inclusion $j''$ with $y = xj''$.
		\end{itemize}
		Then the factorization $f = xj$ is unique.
	\end{prop}
	
	A morphism $f$ in a category with factorization is said to have \textit{image} if $f$ satisfies the condition $(Im)$ of the proposition above. In this case the unique canonical factorization $f = xj$ with the property.,  $(Im)$ is denoted by $ f = f^0j_f$ where $f^0$ is the epimorphic component of $f$ and $j_f$ is inclusion. The unique vertex 
	$Im f = cod f^0 = dom j_f$ is called the image of $f$. Since categories $Set, Grp, $etc., has unique factorization, morphisms in these categories have images. Though the category $Top$ does not have unique factorization, every morphism in $Top$ has image.
	
	\begin{defn}
		Let $F :\mathcal{C} \rightarrow \mathcal{D}$ be a functor. A universal arrow from $d \in \nu \mathcal{D}$ to the functor $F$ is a pair $(c,g)$ where $c \in \nu C$ and $g \in \mathcal{D}(d,F(c))$ such that given any pair $(c',g')$ with $c' \in \nu \mathcal{C}$ and $g \in \mathcal{D}(d,F(c'))$, there is a unique $f \in C(c,c')$ such that $g' = g \circ F (f )$. In this case, the morphism $g$ is universal from $d$ to $F$. A universal arrow from $F$ to $d$ is defined dually.
	\end{defn}
	\begin{exam}
		$Vect_{K}$ be the category of all vector spaces over a field $K$, with arrows linear transformations, then 
		$U : Vect_{K}\rightarrow Set$ sending each vector space $V$ to the set of its elements is a functor and is called the forgetfull functor. For any set $X$ there is a familiar vector space $V_{X}$ with $X$ as the set of basis vectors, 
		consists of all formal $K$-linear combinations of the elements of $X$. Then  $V$ which sends each $x \in X$ into the same $x$ regarded as a vector of $V_{X} $ is also a functor. Consider the arrow $j : X \rightarrow U(V_{X})$. For any other vector space $W$, it is a fact that each function $f : X \rightarrow U(W)$ can be extended to a unique linear transformation $f' : V_{X} \rightarrow W$ with $Uf' \circ j = f$. This familiar fact states that $j$ is a universal arrow from $X$ to $U$.
	\end{exam}
	\begin{defn}
		Let $\mathcal{C}$ be a category and $X_{1}, X_{2}\in v\mathcal{C}$. An object $X$ is a product of $X_{1}$ and $X_{2}$, denoted $ X_{1} \times X_{2}$, if it satisfies this universal property such that:
		there exist morphisms $\pi_{1} : X \rightarrow X_1$, $\pi_2 : X \rightarrow X_2$ such that for every object $Y$ and pair of morphisms $f_1 : Y \rightarrow X_1$, $f_2 : Y \rightarrow X_2$ there exists a unique morphism $f : Y \rightarrow X$ such that the following diagram commutes:
		\begin{center}
			\begin{tikzcd}[column sep=small]
				&Y \arrow[swap]{dl}{f_1} \arrow[dashed]{d}{f} \arrow{rd}{f_2}\\
				X_1 & X_1 \times X_2 \arrow{l}{\pi_{1}} \arrow[swap]{r}{\pi_2} & X_2\\
			\end{tikzcd}
		\end{center}	
		The unique morphism $f$ is called the product of morphisms $f_1$ and $f_2$ and is denoted $< f_1, f_2 >$. The morphisms $\pi_1$ and $\pi_2$ are called the canonical projections.
	\end{defn}
	\begin{exam}
		In the category of sets, the product (in the category theoretic sense) is the cartesian product and in the category of groups, the product is the direct product of groups given by the cartesian product with multiplication defined componentwise.
	\end{exam}
	
	\begin{defn}
		Let $\mathcal{C}$ be a category with zero morphisms. A kernel of a morphism $f: X\to Y$ in $\mathcal{C}$ is an object $K$ together with a morphism $k : K \rightarrow X$ such that 	$f \circ k$ is the zero morphism from $K$ to $Y$.
		Given any morphism $k' : K'\rightarrow  X$ such that $f \circ k'$ is the zero morphism, there is a unique morphism $u : K'\rightarrow  K$ such that $k \circ u = k'$.
	\end{defn}
	
	Kernel of homoprphisms are familiar concept in many categories such as the category of groups or the category of (left) modules over a fixed ring (including vector spaces over a fixed field). To make this explicit, if $f : X \rightarrow Y$ is a homomorphism in one of these categories, and $K$ is its kernel in the usual algebraic sense, then $K$ is a subalgebra of $X$ and the inclusion homomorphism from $K$ to $X$ is a kernel in the categorical sense.
	
	\begin{defn}
		Let $\mathcal{C}$ be a category and $\sigma ^* \mathcal{C}$ be the preorder in $\mathcal{C}$. A cone from $\sigma ^* \mathcal{C}$ to the vertex $d \in \nu \mathcal{C}$ is a map $\gamma: a \mapsto \gamma (a) \in \mathcal{C}(a,d)$ of $\nu \mathcal{C}$ to $\mathcal{C}$ such that whenever $a \subseteq b, j_a^b \gamma(b) = \gamma(a)$.
	\end{defn}
	
	\begin{defn}
		The simplex category $\Delta$ has objects the set of nonempty finite linearly ordered sets $\{[n] | n \geq 0 \}$ and morphisms from $[m]$ to $[n]$ are given by weakly monotone maps. A map $f:\{0, 1, \cdots ,m\} \rightarrow \{0, 1, \cdots ,n\}$ is called weakly monotone if, for every $0 \leq i \leq i'\leq m$, we have $f(i) \leq f(i')$.
	\end{defn}
	The augmented simplex category, denoted by $\Delta_+$ is the category of all finite ordinals and order-preserving maps, thus $ \Delta _{+}=\Delta \cup [-1]$, where$  [-1]=\emptyset $. $[-1]$ is an initial object  and $[0]$ is a terminal object in this category.
	\begin{defn}
		A chain complex $(C, \partial)$ is a graded Abelian group $C = \{C_n\}$ together with an endomorphism $\partial = \{\partial_n\}$ of degree -1, called boundary homomorphism $\partial = \{\partial_n : C_n \rightarrow C_{n-1}\}$, such that $\partial ^2 = 0$. (cf. \cite{david})
	\end{defn}
	
	\section{Category of Chian Bundles}
	We introduced the category of chain bundles in the paper entitled " Category of Chain Bundles" (cf. \cite{PG}). However the following definition is a slightly generalized one.

	\begin{defn}
		Let $\mathcal{C}$ be category with zero. A chain bundle $c$ in the category $\mathcal{C}$ is a sequence \\
		\[\cdots M_{3} \stackrel{S_3}{\rightarrow}M_2 \stackrel{S_2}{\rightarrow} M_1 \stackrel{S_1}{\rightarrow} M_0 = \textbf{0} \]
		where $M_i \in \nu \mathcal{C} $ and $S_i$ be any subset of $Hom(M_{i+1},M_i) $ for all $i$, which also include homsets of the form $Hom(M_{i},M_i)$ and all possible composite of morphisms.
		Let $d$ be the chain bundle \\
		\[\cdots N_{3} \stackrel{T_3}{\rightarrow}N_2 \stackrel{T_2}{\rightarrow} N_1 \stackrel{T_1}{\rightarrow} N_0 = \textbf{0} \]
		A morphism $m: c\to d$ of chain bundles the sequence $m=(f_i,x_i,y_i)$ where $f_i : M_i\to N_i,\,x_i\in S_i$ and $y_i\in T_i$ 
		be such that diagram commutes 
		\begin{center}
			\begin{tikzcd}
				\cdots \arrow[" "]{r} & M_3 \arrow[" S_3"]{r} \arrow{d}[swap]{f_3}
				& M_2 \arrow[" S_2"]{r} \arrow{d}{f_2} & M_1 \arrow{d}{f_1} \arrow["S_1 "]{r} & \textbf{0} \arrow{d}{f_0}\\
				\cdots \arrow[" "]{r}& N_3 \arrow ["T_3 "]{r}  & N_2 \arrow["T_2 "]{r} &  N_1 \arrow["T_1 "]{r} & \textbf{0}  
			\end{tikzcd}
		\end{center}
		that is  $x_i \circ f_{i-1} = f_i\circ y_i$. 
	\end{defn}
	
	Generally a bundle  is a triple $(E, p, B)$ where $E, B$ are sets and $p:E\rightarrow B$ is a map. Each $M_{i+1} \stackrel{S_i}{\rightarrow}  M_i$ in chain bundle can be regarded as a collection of bundles and hence $\cdots M_{3} \stackrel{S_3}{\rightarrow}M_2 \stackrel{S_2}{\rightarrow} M_1 \stackrel{S_1}{\rightarrow} M_0 = \textbf{0} $ is called a chain bundle.\\
	
	Let $\mathcal{C}$ be a category with zero. A category whose objects are chain bundles in $\mathcal{C}$ and morphism are chain bundles maps is the category of chain bundles which we denote by $\mathfrak{CB}_{\mathcal{C}}$. 
	Note that all objects and morphisms in a chain bundle category $\mathfrak{CB}_{\mathcal{C}}$ are objects and morphisms in $\mathcal{C}$  as well and so 
	$\mathfrak{CB}_{\mathcal{C}}$ may regard as a subcategory of $\mathcal{C}$ and further it is easy to observ that given any category $\mathcal{C}$ one can always obtain chain bundle categories $\mathfrak{CB}_{\mathcal{C}}$ of $\mathcal{C}$. The cochain bundles and the category of cochain bundles may be described dually.
	
	Next we proceed to describe the subobject relation and factorization of morphisms in the category $\mathfrak{CB}_{\mathcal{C}}$, whenever the category $\mathcal{C}$ admits subobjects and factorization of morphisms.
	
	\begin{defn}
		Consider the  chain bundle
		$c$: \[\cdots M_{3} \stackrel{S_3}{\rightarrow}M_2 \stackrel{S_2}{\rightarrow} M_1 \stackrel{S_1}{\rightarrow} M_0 = \textbf{0} \in \mathfrak{CB}_{\mathcal{C}}\]
		then, the chain bundle $c'$:\[\cdots M_{3}'\stackrel{S_3'} {\rightarrow}M_2' \stackrel{S_2'}{\rightarrow} M_1' \stackrel{S_1'}{\rightarrow} M_0' = \textbf{0} \]
		with $M_i'$ a subobject of $M_i$ and for each $x' \in S_i' $  there is an  $x \in S_i$ such that $ (j_{M_{i'}}^{M_{i}}x)^0 = x_i'$ is a subchain bundle of $c$. 
	\end{defn}	
	
	Consider chain bundles $c :\cdots \rightarrow M_i \stackrel{S_i}{\rightarrow} M_j \rightarrow \cdots ,\quad  d  :\cdots \rightarrow N_i \stackrel{T_i}{\rightarrow} N_j \rightarrow \cdots \in  \nu \mathfrak{CB}_\mathcal{C}$ and morphism $m: c \rightarrow d$. Since each $f_i : M_i \rightarrow N_i$, admits a factorization of the form $f_i^0j_{N_i'}^{N_i}$ where $N_i' = cod\quad f_i^0$  in $\mathcal{C}$, we obtain a factorization of $m= m^{\circ}J$. For, 
	\begin{center}
		\begin{tikzcd}
			c: \arrow{d} {m} &	\cdots \arrow[" "]{r} & M_i \arrow["S_i "]{r} \arrow{d}[swap]{f_i}
			& M_j \arrow[" "]{r} \arrow{d}{f_j} &  \cdots  \arrow[" "]{r} & M_1\arrow{d}{f_1} \arrow[" "]{r} & \textbf{0} \arrow{d}{0}\\
			d : &	\cdots \arrow[" "]{r}& N_i \arrow ["T_i "]{r}  & N_j \arrow[" "]{r} & \cdots  \arrow[" "]{r} & N_1 \arrow[" "]{r} & \textbf{0}  \\
		\end{tikzcd}
	\end{center}
	each $g \in S_i$ and $ m(g) \in T_i$, $(j_{N_i'}^{N_i}m(g))^0 \in Hom(N_i',Nj')$, take 
	$$S_i' =\{(j_{N_i'}^{N_i}m(g))^0 : f \in S_i \}.$$ 
	
	Define $ m^0 : g \mapsto  (j_{N_i'}^{N_i}m(g))^0 = m(g)' $, then $m(g) = J((j_{N_i'}^{N_i}m(g))^0) =J(m^0(g))$ and so $m = m^0J $ is a factorization. i.e., the category $\mathfrak{CB}_\mathcal{C}$ admits a  factorization as below.
	
	\begin{center}
		\begin{tikzcd}
			c: \arrow{d} {m^0} &	\cdots \arrow[" "]{r} & M_i \arrow["g "]{r} \arrow{d}[swap]{f_i^0}
			& M_j \arrow[" "]{r} \arrow{d}{f_j^0} &  \cdots  \arrow[" "]{r} & M_1\arrow{d}{f_1^0} \arrow[" "]{r} & \textbf{0} \arrow{d}{0}\\
			c' : \arrow{d}{J} & \cdots \arrow[" "]{r} &  N_i' \arrow["m(g)' "]{r} \arrow{d}[swap]{j_i}
			& N_j' \arrow[" "]{r} \arrow{d}{j_j} &  \cdots  \arrow[" "]{r} & N_1' \arrow{d}{j_1} \arrow[" "]{r} & \textbf{0} \arrow{d}{0}\\
			d : &	\cdots \arrow[" "]{r}& N_i \arrow ["m(g) "]{r}  & N_j \arrow[" "]{r} & \cdots  \arrow[" "]{r} & N_1 \arrow[" "]{r} & \textbf{0}  \\
		\end{tikzcd}\\
	\end{center}
	\begin{exam}
		Consider the following chain bundle $\mathfrak{CB_\mathcal{C}}$ where $\mathcal{C}$ is the category $\mathbb{Z}_n$ of finite abelian groups. For $c,d \in \mathfrak{CB_\mathcal{C}}$ the morphism $m: c\to d$
		
		\begin{center}
			
			\begin{tikzcd}
				c: \arrow{d} {m} & \mathbb{Z}_3  \arrow["2a"]{r} \arrow{d}{4} & \mathbb{Z}_6 \arrow["b"]{r} \arrow{d}{4} & \mathbb{Z}_2 \arrow["0"]{r} \arrow{d}{1} & \textbf{0} \arrow{d}{0}\\
				d: & \mathbb{Z}_4 \arrow["2a"]{r} & \mathbb{Z}_8 \arrow["b"]{r} & \mathbb{Z}_2 \arrow["0"]{r} & \textbf{0}\\
			\end{tikzcd}
		\end{center}
		admits the following factorization:
		\begin{center}
			\begin{tikzcd}
				c: \arrow{d} {m^0} & \mathbb{Z}_3  \arrow["2a"]{r} \arrow{d}{4} & \mathbb{Z}_6 \arrow["b"]{r} \arrow{d}{4} & \mathbb{Z}_2 \arrow["0"]{r} \arrow{d}{1} & \textbf{0} \arrow{d}{0}\\
				d': \arrow{d}{J} & \{0\} \arrow["2a"]{r} \arrow{d}{j} & \{0,4\} \arrow["b"]{r} \arrow{d}{j} & \mathbb{Z}_2 \arrow["0"]{r} \arrow{d}{j} & \textbf{0} \arrow{d}{j}\\
				d: & \mathbb{Z}_4 \arrow["2a"]{r} & \mathbb{Z}_8 \arrow["b"]{r} & \mathbb{Z}_2 \arrow["0"]{r} & \textbf{0}\\
			\end{tikzcd}
		\end{center}
	\end{exam}
	
	\subsection{Functors between chain bundle categories}
	\paragraph{}
	Let $\mathcal{C}$ and  $\mathcal{C'}$ be two categories with zero object, then $\mathfrak{CB}_\mathcal{C}$ and $\mathfrak{CB}_\mathcal{C'}$ denote the chain bundle categories. A functor \={F} :  $\mathfrak{CB}_\mathcal{C} \rightarrow \mathfrak{CB}_\mathcal{C'}$ consists of two related functions: the object function \={F} which assigns to each chain bundle $c$ of  $\mathfrak{CB}_\mathcal{C}$ a chain bundle \={F}$c$ of $\mathfrak{CB}_\mathcal{C'}$ and the morphism function which assigns to each morphism $F : c \rightarrow c'$ of  $\mathfrak{CB}_\mathcal{C}$ a morphism \={F}$(F) : $\={F}$c \rightarrow$\={F}$c'$ of $\mathfrak{CB}_\mathcal{C'}$, in such a way that  \={F}$(1_c) = 1_{\text{\={F}c}}$, \={F}$(F \circ G) = \text{\={F}}F \circ \text{\={F}}G$.
	
	\begin{exam}
		Consider the categories $\mathcal{C} = \textbf{Grp}$ and $\mathcal{D} = \textbf{Set}_*$. A functor from $\mathfrak{CB}_\mathcal{C}$ to 
		$\mathfrak{CB}_\mathcal{D}$ defined as follows :\\
		a chain bundle $c:$ \[\cdots G_{3} \stackrel{S_3}{\rightarrow}G_2 \stackrel{S_2}{\rightarrow} G_1 \stackrel{S_1}{\rightarrow} G_0 = \textbf{0} \] in $\mathfrak{CB}_\mathcal{C}$ is mapped to chain bundle $d:$ \[\cdots (G_{3},e_{G_{3}}) \stackrel{S_3}{\rightarrow}(G_{2},e_{G_{2}}) \stackrel{S_2}{\rightarrow} (G_{1},e_{G_{1}}) \stackrel{S_1}{\rightarrow} (G_{0},e_{G_{0}}) = \textbf{0} \] in $\mathfrak{CB}_\mathcal{D}$, where each $G_i$ in $c$ is mapped to corresponding pointed set $(G_{i},e_{G_{i}})$ in $d$ ($e_{G_{i}} $ is the identity element in group  $G_i$) and homomorphisms in $S_i$ of $c$ are mapped to underlying set map is the forgetful functor from $\mathfrak{CB}_\mathcal{C}$ to 
		$\mathfrak{CB}_\mathcal{D}$.
	\end{exam}
	
	\newpage
	\subsection{Chain Bundle Category from Directed Graph and Augmented Simplex Category}
	\paragraph{}
	
	$G$ be a directed graph (with loop at each vertex), then $G$ is regarded as a category with vertices as objects and path between two vertices as morphism. For a vertex $v$ in $G$, the set of all paths ending at $v$ is denoted as $c_v$ and any subset of $c_v$ is a chain bundle in $G$ with end vertex $v$.
	Let $v$ and $w$ be two vertices in $G$ and $c_1, c_2$ be chain bundles in $G$ such that $c_1 \subseteq c_v$ and $c_2 \subseteq c_w$. A morphism from $c_1$ to $c_2$ is the set $P$ of all paths from $v$ to $w$ such that for each path $p \in P$, $c_1^ip \in c_2, \forall c_1^i \in c_1$.
	Thus we can obtain a category of chain bundles from a gaph $G$ by choosing chain bundles in $G$ as objects and chain bundle maps as morphisms and is denoted by $\mathfrak{CB_{\textit{G}}}$. Also this is a category with subobjects where usual inclusion of sets is the subobject relation.\\
	\begin{exam}
		Consider the graph $G$ given below:\\
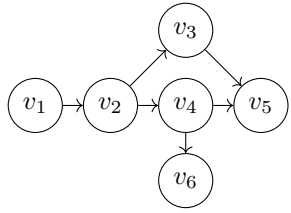
		\begin{center}
			\begin{tikzpicture}[main/.style = {draw, circle}] 
				\node[main] (1) {$v_1$}; 
				\node[main] (2) [right of=1] {$v_2$};
				\node[main] (3) [right of=2] {$v_4$}; 
				\node[main] (4) [above  of=3] {$v_3$};
				\node[main] (5) [right of=3] {$v_5$}; 
				\node[main] (6) [below of=3] {$v_6$};
				\draw[->] (1) -- (2);
				\draw[->] (2) -- (3);
				\draw[->] (2) -- (4);
				\draw[->] (3) -- (5);
				\draw[->] (4) -- (5);
				\draw[->] (3) -- (6);
				
			\end{tikzpicture} 
		\end{center}
		$c_{v_5}=\{v_1v_2v_4v_5,v_1v_2v_3v_5, v_2v_4v_5,v_4v_5,v_3v_5,v_5\}$, $c_{v_2}= \{v_1v_2,v_2\}$, $\{v_1v_2v_3,v_3\}$, $\{v_1v_2v_4v_6, v_2v_4v_6,v_4v_6,v_6\}, c_{v_1}=\{v_1\}$ are chain bundles in $G$. Then  
		$f_1:c_{v_2} \rightarrow  c_{v_5}$ given by
		
		$$v_1v_2 \mapsto v_1v_2v_4v_5,\, v_2 \mapsto v_2v_4v_5$$
		and $f_2:c_{v_2} \rightarrow  c_{v_5}$ given by
		$$v_1v_2 \mapsto v_1v_2v_3v_5, \,v_2 \mapsto v_2v_3v_5$$
		are two chain maps from  $c_{v_2}$ to $  c_{v_5}$.
	\end{exam}

	\paragraph{}
	Consider the augmented simplex category $\Delta_+$. A chain bundle in $\Delta_+$ is of the form $$ \cdots [m_3] \stackrel{S_3}{\rightarrow}  [m_2] \stackrel{S_2}{\rightarrow} [m_1] \rightarrow [0]$$ where $[m_i] $ are ordinals and $S_i$ is any subset of $Hom([m_i],[m_{i-1}]) $ such that $\cdots \leq [m_3]\leq [m_2] \leq [m_1] $.\\
	The chain bundle map between two chain bundles can be defined in a similar way as in Definition \ref{Chainbundle map}. The chain bundle category thus obtained from $\Delta_+$ is denoted as $\mathfrak{CB}_{\Delta_+}$
	
	\section{Categorical Properties of chain bundles}
	In general categorical properties like product, coproduct, kernel, cokernel etc. need not carry over to the bundle categories. However, under certain restrictions some of these categorical properties may carry over to bundle categories.
	
	\vspace{0.2cm}
	Let $\mathcal{C}$ be a category with products. Then product in $\mathfrak{CB}_{\mathcal{C}}$  may be defined as follows,\\
	for the chain bundles $c$ and $d$:\\
	\[c:\cdots M_{3}' \stackrel{S_3'}{\rightarrow}M_2' \stackrel{S_2'}{\rightarrow} M_1'\stackrel{S_1'} {\rightarrow} M_0' = \textbf{0} \]
	\[d: \cdots M_{3} \stackrel{S_3}{\rightarrow}M_2 \stackrel{S_2}{\rightarrow} M_1 \stackrel{S_1}{\rightarrow} M_0 = \textbf{0} \] 
	the product $c \times d $ is the chain bundle  
	\[c \times d:\cdots M_{3}' \times M_{3}  \stackrel{S_3'\times S_3}{\longrightarrow}M_2'\times M_{2}\stackrel{S_2'\times S_2} {\longrightarrow} M_1'\times M_{1} \stackrel{S_1'\times S_1} {\longrightarrow} M_0'\times M_{0}  = \textbf{0} \times \textbf{0}  \]
	For any $F: l \rightarrow c$ and $G: l \rightarrow d $ where $l$ is the chain bundle $$l:\cdots N_{3} \stackrel{T_3}{\rightarrow}N_2 \stackrel{T_2}{\rightarrow} N_1 \stackrel{T_1}{\rightarrow} N_0 = \textbf{0},  $$ there exists a chain bundle map $L : l \rightarrow c \times d $ such that for any $k \in T_i $ and $L(k) \in S_i' \times S_i $ there corresponds  $F(k) \in S_i'$ and $G(k) \in S_i$ such that $L(k) = (F(k),G(k))$.\\
	In a similar manner one can define coproducts in $\mathfrak{CB}_{\mathcal{C}}$.
	\begin{exam}
		$\mathcal{C} $ be the category of submodules of $\mathbb{Z}$ and $ \mathfrak{CB}_{\mathcal{C}} $ be category of chain bundles in $\mathcal{C}$.
		Consider the two chain bundles $$ c:  3\mathbb{Z} \stackrel{\frac{2}{3} a} {\rightarrow} 2\mathbb{Z} \stackrel{\frac{5}{2} b}{\rightarrow} 5\mathbb{Z}\stackrel{0}{\rightarrow} \textbf{0} \quad \text{and} \quad d: 6\mathbb{Z}\stackrel{\frac{2}{3} a'}{\rightarrow} 4\mathbb{Z} \stackrel{\frac{1}{4} b'}{  \rightarrow} \mathbb{Z} \stackrel{0}{\rightarrow} \textbf{0}$$ 
		The product $c \times d$ is the chain bundle 
		$$ c \times d:  3\mathbb{Z} \times 6\mathbb{Z} \stackrel{\frac{2}{3} a,\frac{2}{3} a'} {\rightarrow} 2\mathbb{Z} \times 4\mathbb{Z}   \stackrel{\frac{5}{2} b,\frac{1}{4} b' }{\rightarrow} 5\mathbb{Z} \times\mathbb{Z} \stackrel{0,0}{\rightarrow} \textbf{0} \times \textbf{0}$$
		
		For any $l:  m\mathbb{Z} \stackrel{\frac{n}{m} a''} {\rightarrow} n\mathbb{Z} \stackrel{\frac{p}{n} b''}{\rightarrow} p\mathbb{Z}\stackrel{0}{\rightarrow} \textbf{0}$ in $\nu \mathfrak{CB}_{\mathcal{C}}$ and for any $F: l \rightarrow c$ and $G: l \rightarrow d$ there exists a chain bundle map $L: l \rightarrow c \times d $ such that the following diagram commutes
		\begin{center}
			\begin{tikzcd}
				& l\arrow{dl}[swap]{F} \arrow{d}{L} \arrow{dr}{G}&\\
				c  & c \times d \arrow{l}{\pi_1} \arrow{r}[swap]{\pi_2} & d\\
			\end{tikzcd}
		\end{center}
	\end{exam}
	\vspace{0.2cm}
	Suppose $\mathcal{C}$ be a category having kernels. For $c, d \in \nu \mathfrak{CB_\mathcal{C}}$, let $F : c \rightarrow d$ be a morphism in $\mathfrak{CB_\mathcal{C}}$ whose vertex mapping is $\{f_i\}$. Then kernel of $F$ is a morphism $K : a \rightarrow c$ where $a$ is a chain bundle in $\mathfrak{CB_\mathcal{C}}$ whose vertex mappings are $\{ker(f_i)\}$ and morphism map of $K$ maps each morphism in $T_i$ to the zero morphism in $S_i'$.\\
	\begin{center}
		\begin{tikzcd}
			\arrow{d}{K} a  :  \cdots \arrow[" "]{r} & A_i \arrow[" T_i"]{r} \arrow{d}[swap]{ker{f_i}}
			& A_j \arrow[" "]{r} \arrow{d}{ker{f_j}} &  \cdots \arrow{r} & \textbf{0} \arrow{d}{ker{f_0}}\\
			\arrow{d}{F} c  :  \cdots \arrow[" "]{r}& M_i \arrow ["S_i'"]{r}\arrow{d}[swap]{f_i}  & M_j \arrow[" "]{r}\arrow{d}[swap]{f_j} &  \cdots\arrow{r} & \textbf{0} \arrow{d}{f_0} \\
			d  :  \cdots \arrow[" "]{r} & D_i \arrow[" S_i"]{r}  & D_j \arrow[" "]{r} &  \cdots \arrow{r} & \textbf{0} \\
		\end{tikzcd}\\
	\end{center}  
	Similary one can define cokernels in $\mathfrak{CB_\mathcal{C}}$.
	\begin{exam}
		Consider the chain bundles $ c: \mathbb{Z}_4 \stackrel{a}{\rightarrow} \mathbb{Z}_2 \rightarrow \textbf{0}$, $d : \mathbb{Z}_{12} \stackrel{a}{\rightarrow} \mathbb{Z}_6 \rightarrow \textbf{0}$ and the chain bundle map $F : c \rightarrow d $ given by \\
		\begin{center}
			\begin{tikzcd}
				c: \arrow{d}{F}& \mathbb{Z}_4 \arrow["a"]{r} \arrow{d}{3} & \mathbb{Z}_2 \arrow[" "]{r} \arrow{d}{3} & \textbf{0}\arrow{d}{0}\\
				d :& \mathbb{Z}_{12} \arrow["a"]{r} &  \mathbb{Z}_6\arrow[" "]{r}  & \textbf{0}
			\end{tikzcd}
		\end{center}
		then kernel of $F$ is the following:\\
		\begin{center}
			\begin{tikzcd}
				a: \arrow{d}{Ker F} &	\mathbb{Z}_8 \arrow["a"]{r} \arrow{d}{4} & \mathbb{Z}_4 \arrow[" "]{r} \arrow{d}{2} & \textbf{0}\arrow{d}{0}\\
				c: &	\mathbb{Z}_{4} \arrow["a"]{r} &  \mathbb{Z}_2\arrow[" "]{r}  & \textbf{0}
			\end{tikzcd}
		\end{center}
	\end{exam}
	Let $\mathfrak{CB_\mathcal{C}}$ and $\mathfrak{CB_\mathcal{D}}$ be two chain bundle categories and $\bar{S}$: $\mathfrak{CB_\mathcal{D}} \rightarrow \mathfrak{CB_\mathcal{C}}$ be a functor. For a chain bundle $c$ in $\mathfrak{CB_\mathcal{C}}$ , a universal arrow from $c$ to $\bar{S}$ is a pair $<r,G>$; $r$ is a chain bundle  in $\mathfrak{CB_\mathcal{D}}$ and $G : c \rightarrow \bar{S}r$ such that for every pair $<d,F>$; $ d \in \nu \mathfrak{CB_\mathcal{D}}$ and $F : c \rightarrow \bar{S}d$ there exists a unique arrow $F^\prime : r \rightarrow d$ of $\mathfrak{CB_\mathcal{D}}$ with $G \circ \bar{S}F^\prime = F$.\\
	Consider the category of $Vect_K$ of all vector spaces over a field $K$ and the category $Set$ of sets. $U: Vect_K \rightarrow Set$ be the forgetful functor. For $X \in \nu Set$, universal arrow from $X$ to $U$ is the pair $(V_X,j)$ where $V_X$ is the vector space generated by $X$ and $j: X \rightarrow U(V_X)$ is the inclusion morphism.
	\begin{exam}
		Consider forgetful functor $\bar{U} : \mathfrak{CB_\mathcal{D}} \rightarrow \mathfrak{CB_\mathcal{C}}$ where $\mathcal{C}= Vect_K$ and $\mathcal{D}= Set_*$.
		Let $c: \cdots \rightarrow (X,x) \rightarrow (Y,y) \rightarrow \cdots \rightarrow {\textbf{0}} \in \nu \mathfrak{CB_\mathcal{C}}$. A universal arrow from $c$ to $\bar{U}$ is the pair $<r, \bar{J}>$ where 
		$r: \cdots \rightarrow V_X \rightarrow V_Y \rightarrow {\textbf{0}}$ and $\bar{J}: c \rightarrow \bar{U}(r) $ is defined as follows:\\
		\begin{center}
			\begin{tikzcd}[ampersand replacement=\&]
				c \arrow{d}{\bar{J}}  \&	\cdots \arrow[" "]{r} \& (X,x) \arrow[" "]{r} \arrow{d}{j_X} \& (Y,y) \arrow[" "]{r} \arrow{d}{j_Y} \& \textbf{0}\arrow{d}{0}\\
				\bar{U}(r) \&	\cdots \arrow[" "]{r} \& (\bar{U}(V_X,0_{V_X})) \arrow[" "]{r} \&  (\bar{U}(V_Y,0_{V_Y})) \arrow[" "]{r} \& \textbf{0}
			\end{tikzcd}
		\end{center}
		where $j_X: (X,x) \rightarrow \bar{U}(V_X,0_{V_X})$ is defined by $j_X(x)= 0, j_X(X-x) = I(X-x) $ and for $f: X \rightarrow Y $ with$ f(x)= y$, $\bar{J}(f) : \bar{U}(V_X) \rightarrow \bar{U}(V_Y)$ is given by\\
		$$\bar{J}(f)(j_X(x)) = j_Y(y)$$\\
		$$\bar{J}(f)(\bar{U}(V_X)-X) = 0_{V_Y}$$\\
		$$\bar{J}(f)(X-j_X(x)) = f(X-j_X(x))$$
	\end{exam}
	Let $S: \mathcal{D} \rightarrow \mathcal{C}$ be a functor and $\bar{S}$ be the functor from $\mathfrak{CB_\mathcal{D}} $ to $ \mathfrak{CB_\mathcal{C}}$ induced by $S$. Suppose that for each $C_i \in \nu \mathcal{C}$ there is a universal arrow $<D_i,g_i>$ form $C_i$ to $S$, where $D_i \in \nu \mathcal{D}$ and $g_i : C_i \rightarrow S(D_i)$ is a retraction. Now, for $c: \cdots \rightarrow C_i \stackrel{A_i}{\rightarrow} C_j \rightarrow \cdots \rightarrow \textbf{0} \in \nu \mathfrak{CB_\mathcal{C}}$, we have  $<d,G>$ where $d: \cdots \rightarrow D_i \stackrel{Hom(D_i,D_j)}{\rightarrow} D_j \rightarrow \cdots \rightarrow \textbf{0} \in \nu \mathfrak{CB_\mathcal{D}}$ and $G: c \rightarrow\bar{S}d$ with $\nu G = \{g_i\}$ and for $x \in A_i$, $G(x) = g_i^{-1}xg_j$ is a universal arrow from $c$ to $\bar{S}$.

	\section{Category of Chains}
	Next we proceed to discuss the chains in the category of chain bundles $\mathfrak{CB_\mathcal{C}}$. A chain is a chain bundle having exactly one morphism in each homset. A category whose objects are chains and morphisms are appropriate chain maps is called category of chains.
	\begin{defn}
		Let $\mathcal{C}$ be a category with zero object. A chain in $\mathcal{C}$ is \[\cdots M_{3} \stackrel{s_3}{\rightarrow}M_2 \stackrel{s_2}{\rightarrow} M_1 \stackrel{s_1}{\rightarrow} M_0 = \textbf{0} \]
		where $M_i \in \nu \mathcal{C} $ and $s_i$ is a morphism in $Hom(M_{i+1},M_i) \quad \forall i$, which also consists of morphisms $1_{M_i}$ and all possible composite of morphisms.
	\end{defn}
	\begin{defn} 
		A chain map between two chains in $\mathcal{C}$ is a functor $F$ between the two whose vertex map $\nu F = \{f_i: M_i \rightarrow N_i\} $ is a sequence of morphisms in $\mathcal{C}$ such that the below diagram commutes. 
		\begin{center}
			\begin{tikzcd}
				\cdots \arrow[" "]{r} & M_3 \arrow[" s_3"]{r} \arrow{d}[swap]{f_3}
				& M_2 \arrow[" s_2"]{r} \arrow{d}{f_2} & M_1 \arrow{d}{f_1} \arrow["s_1 "]{r} & \textbf{0} \arrow{d}{f_0}\\
				\cdots \arrow[" "]{r}& N_3 \arrow ["s_3' "]{r}  & N_2 \arrow["s_2' "]{r} &  N_1 \arrow["s_1' "]{r} & \textbf{0}  \\
			\end{tikzcd}\\
			
		\end{center}
	\end{defn}
	\begin{defn}
		Given a chain $c$ in the category  $\mathcal{C}$ 
		\[c:\cdots M_{3} \stackrel{s_3}{\rightarrow}M_2 \stackrel{s_2}{\rightarrow} M_1 \stackrel{s_1}{\rightarrow} M_0 = \textbf{0} \]
		then the chain \[ c':\cdots M_{3}'\stackrel{s_3'} {\rightarrow}M_2' \stackrel{s_2'}{\rightarrow} M_1' \stackrel{s_1'}{\rightarrow} M_0' = \textbf{0} \]
		with $M_i'$ a subobject of $M_i$ and $(s_{i}')^0 = (j_{M_{i'}}^{M_{i}}s_i)^0$is a subchain of $c$. 
	\end{defn}
	\begin{exam}
		Consider the category $\mathcal{C}$ of graded abelian groups and the category $ \mathfrak{CB}_{\mathcal{C}}$ of chain bundles in $\mathcal{C}$. Define the condition for choosing chain from a chain bundle in $ \mathfrak{CB}_{\mathcal{C}}$ as follows:\\
		Let  $c: \cdots C_{n+1} \rightarrow C_{n} \rightarrow C_{n-1}\rightarrow \cdots \rightarrow \textbf{0} $ is a chain bundle in  $ \mathfrak{CB}_{\mathcal{C}}$ where $\{C_n\}$ are abelian groups, then choose one  $\partial_i  $ from each $Hom (C_{i+1},C_{i})$ such that $\partial_{i+1}\circ \partial_{i} = 0 $. If we can choose such a $\partial_i  $ from each $Hom (C_i,C_{i+1}) \quad \forall i$ we get a chain of the form \\
		$$c_1: \cdots C_{n+1} \stackrel{\partial_{n+1}}{\rightarrow} C_n \stackrel{\partial_n}{\rightarrow} C_{n-1} \stackrel{\partial_{n-1}}{\rightarrow} \cdots \rightarrow \textbf{0}$$\\
		Chains in the chain bundle category together with chain maps is a category which we call category of chains.
	\end{exam}
	Let $\mathcal{C}$ be a category with factorization and $\Gamma $ be a category of chains obtained from $\mathfrak{CB_\mathcal{C}}$.  Consider the chains  $c :\cdots \rightarrow M_i \stackrel{s_i}{\rightarrow} M_j \rightarrow \cdots $ and $ d  :\cdots \rightarrow N_i \stackrel{t_i}{\rightarrow} N_j \rightarrow \cdots \in  \nu \Gamma$, and the chain map $F: c \rightarrow d$ with $\nu F =\{f_i : M_i \rightarrow N_i\}$. Since each $f_i$ admits a factorization in $\mathcal{C}$, $F$ may admit a factorization of the form $F=F^0J$ as in case of factorization of cahin bundle map depending on the choice of category of chains.\\
	\paragraph*{}
	In the following we provide an example of a category of chains with subobjects and factorization.
	\begin{exam}
		Let \textbf{Csim} denote category of simplicial complexes and simplicial maps. Consider the following simplicial complex.
		
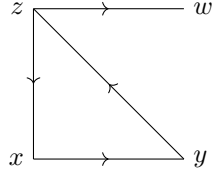
\begin{center}	
			\begin{tikzpicture}[node distance=2cm]
				\coordinate[label=right:$w$] (tw) {};
				\coordinate[left=of tw,label=left:$z$] (tz) {};
				\coordinate[below=of tz,label=left:$x$] (bx) {};
				\coordinate[below=of tw,label=right:$y$] (by) {};
				\begin{scope}[decoration={markings,mark=at position 0.5 with {\arrow{>}}}] 
					\draw[postaction={decorate}] (tz) -- node[auto,swap] {} (tw);
					\draw[postaction={decorate}] (tz) -- node[auto] {} (bx);
					\draw[postaction={decorate}] (bx) -- node[auto,swap] {} (by);
					\draw[postaction={decorate}] (by) -- node[auto] {} (tz);
					
				\end{scope}
				
			\end{tikzpicture}
		\end{center}
		Let $\mathcal{C'}$ be the subcategory of \textbf{Csim} whose objects are set of all subcomplexes of above simplicial complex and morphisms are simplicial maps between them.  Consider the homology functor $H_*(-,\mathbb{Z}): \textbf{Csim} \rightarrow \textbf{Ab}^{\mathbb{Z}}$, where $\textbf{Ab}^{\mathbb{Z}}$ denotes category of graded abelian groups and related morphisms ( see cf.\cite{david}). $Im(j_{\mathcal{C'}}^{\textbf{Csim}}H_*(-,\mathbb{Z}))$ is a  subcategory of $\textbf{Ab}^{\mathbb{Z}}$ denoted by  $\mathcal{C}$. Let $\mathfrak{CB_\mathcal{C}}$ be category of chain bundles in $\mathcal{C}$. By chosing chain complexes from $\mathfrak{CB_\mathcal{C}}$  we obtain a category $\Gamma$ with subobjects and factorization. Being a (chain) subcomplex is a subobject relation and any chain map $F$,
		
		\begin{center}
			\begin{tikzcd}
				c: \arrow{d}{F}	\cdots \arrow[" "]{r} & C_2' \arrow[" \partial_2"]{r} \arrow{d}[swap]{f_2}
				& C_1' \arrow[" \partial_1"]{r} \arrow{d}{f_1} & C_0' \arrow{d}{f_0} \arrow["\partial_0 "]{r} & \textbf{0} \arrow{d}{0}\\
				d:	\cdots \arrow[" "]{r}& C_2 \arrow ["\partial_2' "]{r}  & C_1 \arrow["\partial_1' "]{r} &  C_0 \arrow["\partial_0' "]{r} & \textbf{0}  \\
			\end{tikzcd}\\
			
		\end{center}
		admits a factorization given by
		\begin{center}
			\begin{tikzcd}
				c: \arrow{d}{F^0}	\cdots \arrow[" "]{r} & C_2' \arrow[" \partial_2"]{r} \arrow{d}[swap]{f_2^0}
				& C_1' \arrow[" \partial_1"]{r} \arrow{d}{f_1^0} & C_0' \arrow{d}{f_0^0} \arrow["\partial_0 "]{r} & \textbf{0} \arrow{d}{0}\\
				c': \arrow{d}{J}	\cdots \arrow[" "]{r} & F^0(C_2') \arrow[" \partial_{2_{|C_2}}'"]{r} \arrow{d}[swap]{j}
				& F^0(C_1') \arrow[" \partial_{1_{|C_1}}'"]{r} \arrow{d}{j} & F^0(C_0') \arrow{d}{j} \arrow["\partial_{0_{|C_0}}' "]{r} & \textbf{0} \arrow{d}{0}\\
				d:	\cdots \arrow[" "]{r}& C_2 \arrow ["\partial_2' "]{r}  & C_1 \arrow["\partial_1' "]{r} &  C_0 \arrow["\partial_0' "]{r} & \textbf{0}  \\
			\end{tikzcd}\\
			
		\end{center}
		In particular, consider $K_1, K_2 \in \nu \mathcal{C'}$ where \\
		$$K_1 = \{\{x,y,z,w\},\{(xyzw),(xyz),(xy),(yz),(zx),(zw),(x),(y),(z),(w)\}\}$$
		and $$K_2 = \{\{x,y,z,w\},\{(yz),(zx),(zw),(x),(y),(z),(w)\}\}$$\\
		Let $f: K_1 \rightarrow K_2$ be the simplicial map given by\\
		$$ 
		x \mapsto y, 
		y \mapsto z, 
		z \mapsto z, 
		w \mapsto z $$
		$C(K_1) = \{C_n(K_1)\}$ and $ C(K_2) = \{C_n(K_2)\}$ be chain complexes corresponding to $K_1$ and $K_2$. $C_n(f): C_n(K_1) \rightarrow C_n(K_2)$ is given by
		\begin{equation*}
			C_n(f)(\{x,y,z,w\}) =\begin{cases}
				\{f(x),f(y),f(z),f(w)\}, & \text{if $f(x)$ are distinct for different $x$}.\\
				0, & \text{otherwise}.
			\end{cases}
		\end{equation*}
		
		Consider the chain map
		$C(f) : C(K_1) \rightarrow C(K_2)$  in $\Gamma$ induced by $f$:
		\begin{center}
			\begin{tikzcd}
				C(K_1): \arrow{d}{C(f)} & 0 \arrow["0 "]{r} \arrow{d}[swap]{0} & <xyz> \arrow[" \partial_3'"]{r} \arrow{d}[swap]{0} & <xy,yz,zx,zw> \arrow[" \partial_2'"]{r} \arrow{d}[swap]{C_2(f)}
				& <x,y,z,w> \arrow[" \partial_1'"]{r} \arrow{d}{C_1(f)} & \textbf{0} \arrow{d}{0}\\
				C(K_2):	& 0 \arrow["0 "]{r} & 0 \arrow[" \partial_3"]{r}  & <yz,zx,zw> \arrow[" \partial_2"]{r} 
				& <x,y,z,w> \arrow[" \partial_1"]{r}  & \textbf{0} \\
			\end{tikzcd}\\
			
		\end{center}
		Factorization of $C(f) : C(K_1) \rightarrow C(K_2)$ is given by:
		\begin{center}
			\begin{tikzcd}
				C(K_1): \arrow{d}{C(f)^0} & 0 \arrow["0 "]{r} \arrow{d}[swap]{0} & <xyz> \arrow[" \partial_3'"]{r} \arrow{d}[swap]{0} & <xy,yz,zx,zw> \arrow[" \partial_2'"]{r} \arrow{d}[swap]{C_2(f)^0}
				& <x,y,z,w> \arrow[" \partial_1'"]{r} \arrow{d}{C_1(f)^0} & \textbf{0} \arrow{d}{0}\\
				C(K_3): \arrow{d}{J} & 0 \arrow["0 "]{r} \arrow{d}[swap]{0} & 0 \arrow[" \partial_3''"]{r} \arrow{d}[swap]{0} & <yz> \arrow[" \partial_2''"]{r} \arrow{d}[swap]{j}
				& <y,z> \arrow[" \partial_1''"]{r} \arrow{d}{j} & \textbf{0} \arrow{d}{0}\\
				C(K_2):	& 0 \arrow["0 "]{r} & 0 \arrow[" \partial_3"]{r}  & <yz,zx,zw> \arrow[" \partial_2"]{r} 
				& <x,y,z,w> \arrow[" \partial_1"]{r}  & \textbf{0} \\
			\end{tikzcd}\\
			where $C(K_3)$ is the chain complex corresponding to simplicial complex $K_3 = \{\{y,z\},\{(yz),(y),(z)\}\}$
		\end{center}
	\end{exam}

\end{document}